\documentclass{amsart}
\usepackage{amscd}
\usepackage{wrapfig}
\usepackage{amsmath,amsfonts,amsthm,amscd,amssymb,latexsym,comment}
\usepackage{enumerate,graphicx,psfrag,subfigure}

\usepackage{amsmath,amsfonts,amsthm,amscd,amssymb,latexsym,comment}
\usepackage{enumerate,graphicx,psfrag,subfigure,changebar,oldgerm}

\newtheorem{thm}{Theorem}
\newtheorem{theorem}[thm]{Theorem}
\newtheorem{example}{Example}
\newtheorem{lemma}[thm]{Lemma}
\newtheorem{corollary}[thm]{Corollary}

\newtheorem{definition}{Definition}
\newtheorem{remark}{Remark}
\newtheorem{question}{Question}

\begin{document}
\title{A Gathering Process in Artin Braid
Groups\footnote{Authors were supported by the RFFI grant
N02-01-00192}}

\author{Evgeinj S. Esyp}

\address{Omsk Branch of Institute of Mathematics \\
(Siberian branch of Russian Academy of Science) \\ Pevtsova st. 13,
Omsk, 644099, Russia; \textsf{email:} esyp@iitam.omsk.net.ru }

\author{Ilya V. Kazachkov}

\address{Department of Mathematics and Statistics \\
McGill University \\ 805 Sherbrooke st. West \\
Montreal, QC, H3A 2K6, Canada;
 \textsf{email:} Ilya.Kazachkov@gmail.com }

\maketitle

\begin{abstract}
In this paper we construct a gathering process by the means of which
we obtain new normal forms in braid groups. The new normal forms
generalise Artin-Markoff normal forms and possess an extremely
natural geometric description. In the two last sections of the paper
we discuss the implementation of the introduced gathering process
and the questions that arose in our work. This discussion leads us
to some interesting observations, in particular, we offer a method
of generating a random braid.
\end{abstract}


\tableofcontents

\section{Preliminaries}

Recall, that Artin braid group on $n+1$ strands is the group given
by the following generators and relations:
\begin{equation} \label{eq:B}
B_{n+1}= \left<x_1, \ldots, x_n| \left[ x_i, x_j \right]=1, |i-j|
\ge 2; x_ix_{i+1}x_i=x_{i+1}x_ix_{i+1} \right>.
\end{equation}

\begin{figure}[h]
  \centering
  \includegraphics[keepaspectratio,width=12.6cm, height=2.6cm]{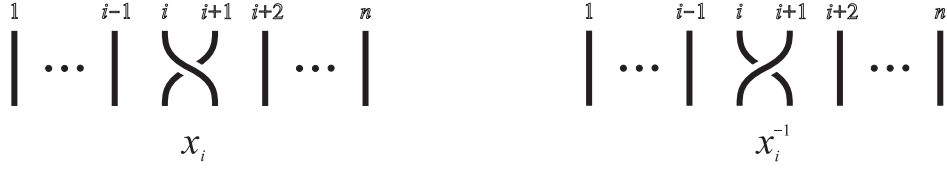}
\caption{The diagrams of generators} \label{pic:cross}
 \end{figure}

The braid group admits a geometric presentation which we shall use
throughout this paper. We refer to \cite{Dehornoy, DehSurv,
Vershinin} for details. The idea is to associate to every word in
the $x_i$'s a plane diagram, which is obtained by successively
concatenating the diagrams of letters, they are given by Figure
\ref{pic:cross}.

\begin{figure}[h]
  \centering
  \includegraphics[height=4cm, width=10cm, keepaspectratio]{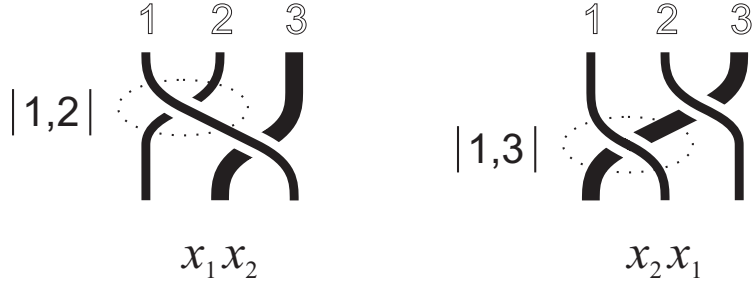}
 \caption{Correspondence of crossings
and letters} \label{pic:1}
 \end{figure}

We shall make use of another definition, closely linked to the
geometric presentation of braids. If two strands cross, then this
we call a crossing (on Figure \ref{pic:cross} the crossings of the
$i$-th and the $i+1$-th strands are shown). To every letter $y$ of
a braid word $\xi \in B_{n+1}$ we associate a `crossing' of two
strands. Consider the diagram of $\xi$. The letters of $\xi$ are
in one-to-one correspondence with the crossings of the diagram of
$\xi$. So, depending on the position of the letter $y$ in the word
$\xi$ the crossing corresponding to $y$ may be different. For
instance, to the letter $x_1$ in the word $\xi_1= x_1x_2$
corresponds the crossing of the first and the second strand, while
in the word $\xi_2=x_2x_1$ the crossing linked to $x_1$ is the
crossing of the first and the third strands (see Figure
\ref{pic:1}). Depending on the sign of the generator the sign of
crossing can be positive or negative.

\begin{definition}[notation of crossing]
Let $\mid r,s \mid^{\epsilon}$ denote the crossing of the $r$-th
and the $s$-th strands. Thus the crossing is completely defined
when we fix the numbers $r$ and $s$ and the sign. We say that a
crossing is negative, whenever the corresponding letter of the
word is the inverse of one of the generators $x_1, \ldots, x_n$
and positive otherwise. Thus, the natural notion of the
\emph{`sign'} $\epsilon=\pm 1$ of a crossing arises.
\end{definition}
From the foregoing discussion we conclude that braid words can be
treated as sequences of crossings.

\begin{figure}[h]
\centering
  \includegraphics[height=4cm, width=10cm, keepaspectratio]{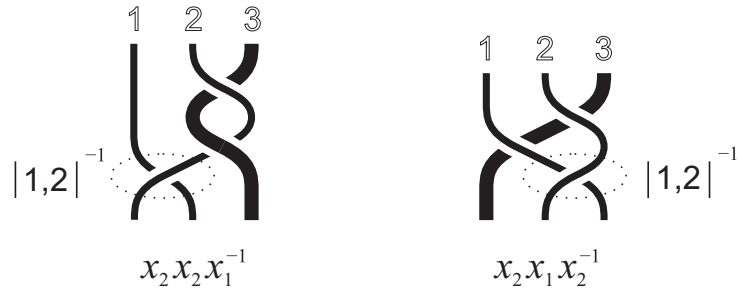}
\caption{Correspondence of crossings and letters} \label{pic:2}
\end{figure}

Note that any crossing, consider $\mid 1, 2\mid^{-1}$, for
instance, can correspond to different generators $x_1, \ldots,
x_n$, depending on the word and the position of a generator in the
word (see Figure \ref{pic:2}).

Furthermore, we offer the reader to check that the following
remark holds

\begin{remark}
Let $\mid i, j \mid ^{\epsilon}$ be an arbitrary crossing. Then
there exists a braid word $w$ and a letter $y$ in $w$ so that the
crossing corresponding to $y$ in $w$ is $\mid i, j \mid
^{\epsilon}$, here $y \in X \cup X^{-1}$.
\end{remark}

We next formulate the main result of the current paper.

\begin{theorem}[Normal form of elements of $B_{n+1}$] \label{thm:NF}
Let $B_{n+1}$ be the braid group on $n+1$ strands. Every element
$w \in B_{n+1}$ can be uniquely written in the form
\begin{equation} \label{eq:NF}
w=x_1^m \cdot w_3(x_1,x_2) \cdot w_4(x_1,x_2,x_3) \cdots
w_{n+1}(x_1, \dots ,x_n);\; m \in \mathbb{Z}
\end{equation}
Where for every $3 \ge k \le n+1$ the words $w_k$  are freely
reduced and the crossings involved into $w_k$ are the crossings of
the form $\mid i,k \mid$ {\rm(}$k=3, \dots, n+1; i<k${\rm)} only.
\end{theorem}

This is roughly saying that for given braid there exists a unique
braid, which is equivalent to the initial one and is `constructed'
in the following way: first one entangles the first and the second
strand $m$ times and the first and the second strands do not cross
further. Then one entangles the third strand with the second and
the first (this is $w_3$) so that furtheron there are no crossings
of the first and the second  strands with the third strand, then
one entangles the fourth strand with the third, the second and the
first (this is $w_4$) and so on.

\begin{center}
\begin{figure}[h]
    \centering
  \includegraphics[height=5cm, width=11cm, keepaspectratio]{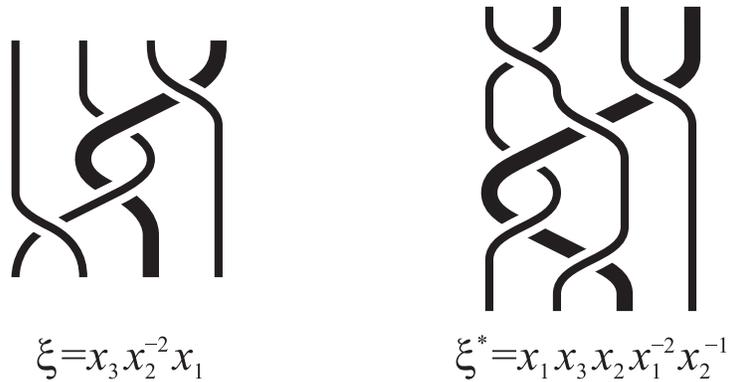}
    \caption{Example 1. Normal form of $\xi = x_3x_2^{-2}x_1$}
\label{pic:expl1}
\end{figure}
\end{center}

\begin{example}
Let $\xi \in B_4$, $\xi = x_3x_2^{-2}x_1$. Then normal form
{\rm(\ref{eq:NF})} of the word $\xi$ is
$\xi^*=x_1x_3x_2x_1^{-2}x_2^{-1}$, see Figure
{\rm\ref{pic:expl1}}.
\end{example}

In the case of $B_3$ Theorem {\rm \ref{thm:NF}} can be
reformulated as follows:

\begin{corollary} \label{cor:B3}
 Every element $w \in B_{3}$ can be
uniquely written in the form
\[
w=x_1^m \cdot (x_2^{k_1}  x_1^{k_2} x_2^{k_3} \cdots x_l^k)
\]
here $l=1,2$, $k, m \in \mathbb{Z}$, $k_r$ is even if and only if
$r$ is even and is odd if and only if $r$ is odd.
\end{corollary}

\begin{remark}
We draw reader's attention to the fact that in the normal form
given in the above Corollary the power of the first letter $x_1$
has an arbitrary degree $m$ and that the terminal letter of $w$
may be any letter {\rm(}$x_1$ or $x_2${\rm)} in any power $k$.

Consider the word $w=x_2^{k_1} x_1^k \in B_3$, where $k_1$ is an
odd integer and $k\in \mathbb{Z}$. This braid is in normal form,
since only the third strand is entangled. If we now look at $w$ as
a word written as in Corollary {\rm\ref{cor:B3}}, we have $m=0$
and, since $x_1$ is the terminal letter of $w$, it may be in an
arbitrary power (not necessarily even).

Moreover, the word $x_1^m x_2^k$, where $m,k \in \mathbb{Z}$, is
in normal form.
\end{remark}

In  paper \cite{Artin} E.~Artin (and A.~Markov in \cite{Markov})
proves that every pure braid can be transformed into the form
(\ref{eq:NF}). He suggests to introduce a normal form in the braid
group $B_{n+1}$ as follows. Every element of $B_{n+1}$ can be
written as a product of a pure braid and a coset representative of
$B_{n+1}$ modulo the group of pure braids $I_{n+1}$. Normal form
(\ref{eq:NF}), therefore, generalises Artin-Markoff normal form
(see \cite{Artin} and \cite{Markov}) to all (not necessarily pure)
braids and admits a natural geometric description. Moreover, in
Section \ref{sec:proof}, we construct a gathering process which
transforms a given word into the form (\ref{eq:NF}).

As mentioned above in the case when $w$ is a pure braid, normal
form  (\ref{eq:NF}) coincides with Artin-Markoff normal form,
given by Theorem 17, \cite{Artin} (see as well \cite{Markov}). For
the sake of convenience and completeness we expose this theorem
below.

By the definition set
\[
A_{i,i+1}=x_i^{2}, \quad i=1,\ldots, n.
\]
For $i<j, i,j=1,\ldots,n+1$ set
\[
A_{j,i}=A_{i,j}=x_i^{-1} \cdots x_{j-1}^{-1} \cdot x_{j}^2 \cdot
x_{j-1}\cdots x_i=x_{j-1} \cdots x_{i+1} \cdot x_i^{2}\cdot
x_{i+1}^{-1} \cdots x_{j-1}^{-1}
\]

\begin{theorem}[\cite{Artin}, Theorem 17] \label{thm:art}
The $A_{i,k}$'s are generators of the group of pure braids. Every
pure braid can be uniquely written in the form:
\begin{equation} \label{eq:art}
A=U_1\cdots U_{n-1}
\end{equation}
where each $U_j$ is a uniquely determined power product of the
$A_{i,j}$ using only those with $i>j$.
\end{theorem}

\begin{center}
\begin{figure}[h]
    \centering
  \includegraphics[height=5cm, width=11cm, keepaspectratio]{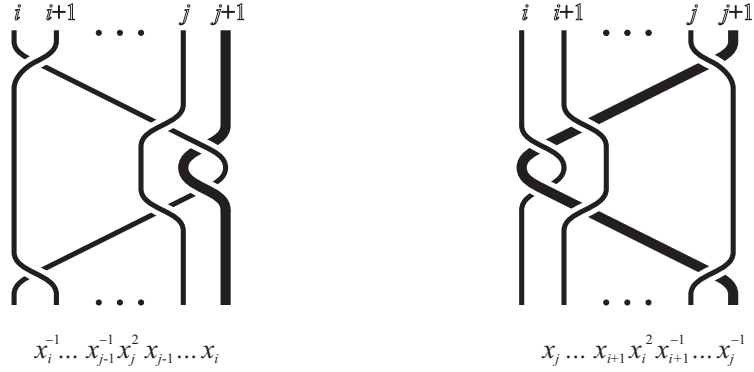}
    \caption{$A_{i,j}=A_{j,i}$}
\label{pic:Aij}
\end{figure}
\end{center}

\begin{remark}
In {\rm \cite{Artin}} the author reads the words from right to
left and draws diagrams moving upwards. The diagrams {\rm(}in our
interpretation{\rm)} of elements $A_{j,i}=A_{i,j}$ take the form
shown on Figure {\rm\ref{pic:Aij}}.

In {\rm \cite{Artin}} the author introduces the elements
$A_{i,j}$'s to make normal form {\rm(\ref{eq:NF})} a unique word
in the $A_{i,j}$'s {\rm(}however the elements $A_{i,j}$'s can be
presented by different braid words{\rm)}. Choose the word
representative of the element $A_{i,j}$ to be $x_{j-1} \cdots
x_{i+1} \cdot x_i^{2}\cdot x_{i+1}^{-1} \cdots x_{j-1}^{-1}$,
where $j>i$. Notice next that in this case the braid word
$A_{i,j}$ entangles the $j$-th strand only, leaving the other $n$
strands unentangled (see Figure \ref{pic:Aij}).

Finally, minding the insignificant difference in the geometric
interpretation of braid words and noticing that the elements
$U_{n-1}, \ldots, U_1$ in {\rm(\ref{eq:art})} are uniquely
determined elements of $B_{n+1}$ {\rm(}in {\rm\cite{Artin}} the
$U_j$'s are words in the symbols $A_{i,j}$'s{\rm)}, we conclude
that for pure braids the elements $x_1^m, w_3, \ldots, w_{n+1}$ in
normal form {\rm(\ref{eq:NF})} coincide with, correspondingly,
$U_{n-1}, \ldots, U_1$ in Equation {\rm(\ref{eq:art})}.
\end{remark}

\section{Proof of Theorem \ref{thm:NF}} \label{sec:proof}

In this section we prove Theorem \ref{thm:NF}. We also provide an
algorithm (a gathering process) for constructing normal form
(\ref{eq:NF}).

We first prove that every braid word $w \in B_{n+1}$ can be taken
to the word in the form (\ref{eq:NF}). And then prove that there
exists only a unique word of the form (\ref{eq:NF}) equal to $w$
in $B_{n+1}$.

We say that a crossing that involves the $n+1$-th strand is
\emph{big}. Otherwise we term a crossing \emph{small}.

In the following by $\xi_1\equiv\xi_2$ we shall mean the equality
in the free group and by $\xi_1 \cdot \xi_2$ we denote
cancellation-free multiplication of two words.

\subsection{Existence} \label{ss:exist}
Below we introduce a gathering process that takes the word $w$
into the form (\ref{eq:NF}). First we transform the word $w$ into
the form in which every big crossing is gathered in the  end of
the word, $w=\xi_n w_{n+1}$. Then we apply a similar process to
$\xi_n$ and take it in the form in which every crossing that
involves the $n$-th strand is collected in the end of the word,
$w=\xi_n w_{n+1}=\xi_{n-1} w_n w_{n+1}$ and so on.

\begin{figure}[p]
\centering
  \includegraphics[height=22cm, width=12.6cm, keepaspectratio]{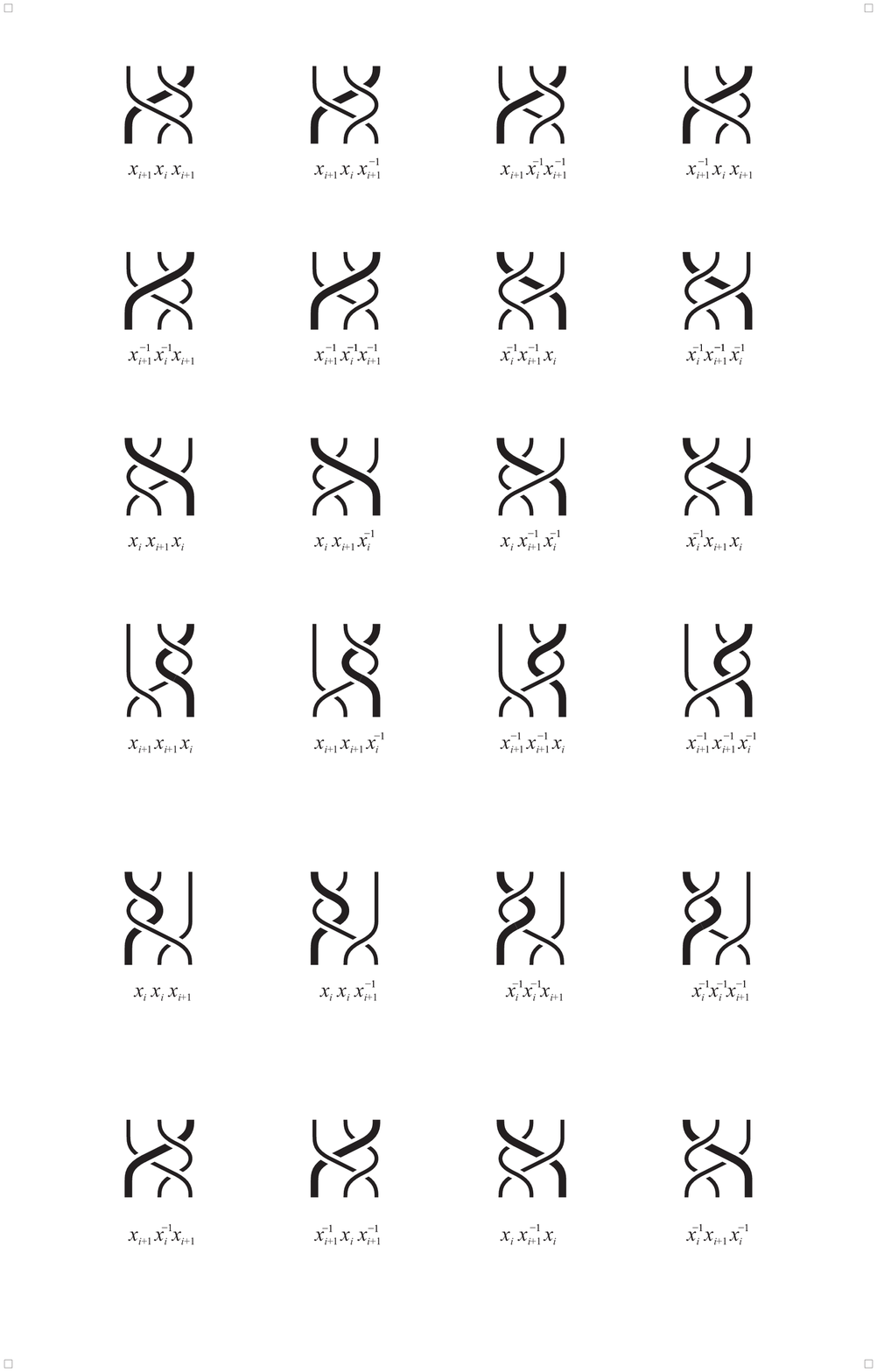} \caption{Possible variants}
    \label{pic:pv}
\end{figure}

\begin{figure}[p]
\centering
  \includegraphics[height=22cm, width=12.6cm, keepaspectratio]{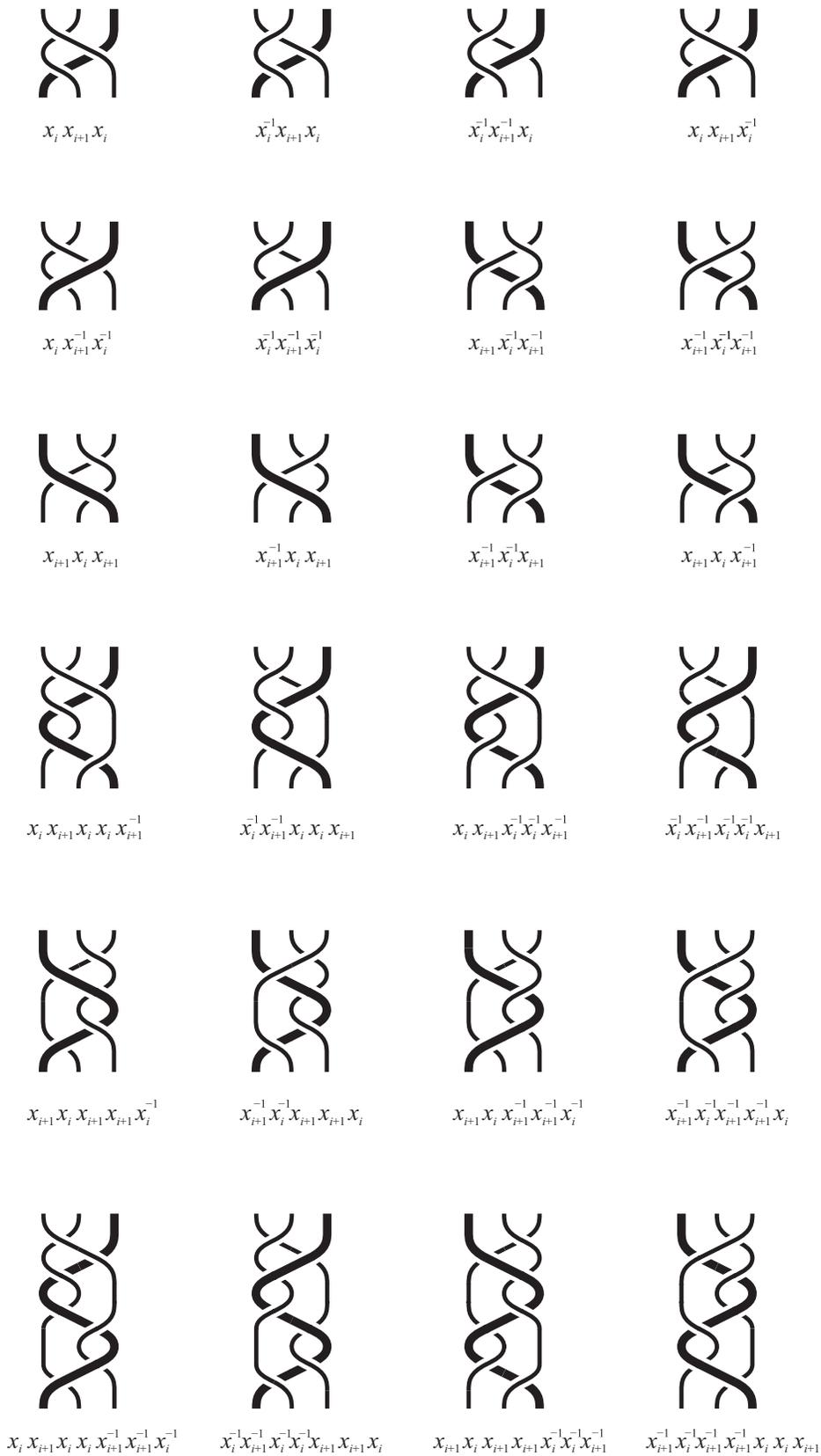} \caption{Replacements}
  \label{pic:repl}
\end{figure}

Before we begin to give a formal description of our algorithm we
give a less formal description of the idea it uses. Given an
arbitrary word $w\in B_n$ we consider it as a sequence of
crossings. Take the first small crossing $\rho$ in $w$ so that
there are big crossings preceding $\rho$ (i.e. the corresponding
letters are closer to the beginning of $w$). The word $w$ then has
the following form $T\cdot U \cdot V \cdot W$. Where in $T$ all
the crossings are small, in $U$ all the crossings are big, $V$ is
a letter and the corresponding crossing is $\rho$ and $W$ is
simply the rest of the word. We then, using transformations given
by Figures \ref{pic:pv} and \ref{pic:repl}, move this crossing
upwards (preserving the braid) till there are no big crossings in
$w$ preceding $\rho$. Applying the same procedure to every small
crossing in $w$ we get its decomposition in the form $w=\xi_n
w_{n+1}$

We begin to construct a gathering process that transforms an
arbitrary braid word $w$ into the form $w=\xi_n w_{n+1}$, where
every crossing in $w_{n+1}$ is big and every crossing in $\xi_n$
is small.

\begin{figure}[h]
\centering
  \includegraphics[height=3cm, width=6cm, keepaspectratio]{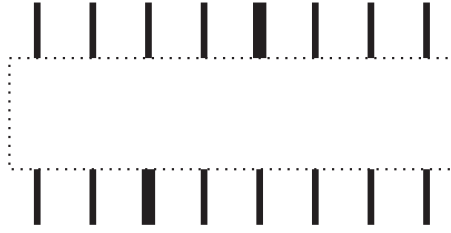} \caption{An arbitrary
  braid}\label{pic:3}
\end{figure}

\begin{figure}[h]
\centering
  \includegraphics[height=3cm, width=6cm, keepaspectratio]{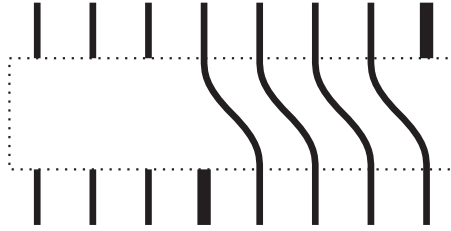} \caption{A braid in which every
crossing is big} \label{pic:4}
\end{figure}

\begin{figure}[h]
\centering
  \includegraphics[height=3cm, width=6cm, keepaspectratio]{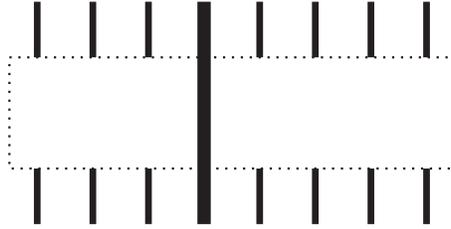}
  \caption{A braid in which every crossing is small}
\label{pic:5}
\end{figure}

By the diagram shown on Figure \ref{pic:3} we denote an arbitrary
braid on $n+1$ strands. Let us agree that the $n+1$-th strand is
designated bold.
By the diagram shown on Figure \ref{pic:4} we denote a braid, in
which no small crossings occur.
And by the diagram shown on Figure \ref{pic:5} we denote a braid,
in which no big crossings occur.

\begin{figure}[h]
\centering
  \includegraphics[height=12cm, width=8cm, keepaspectratio]
{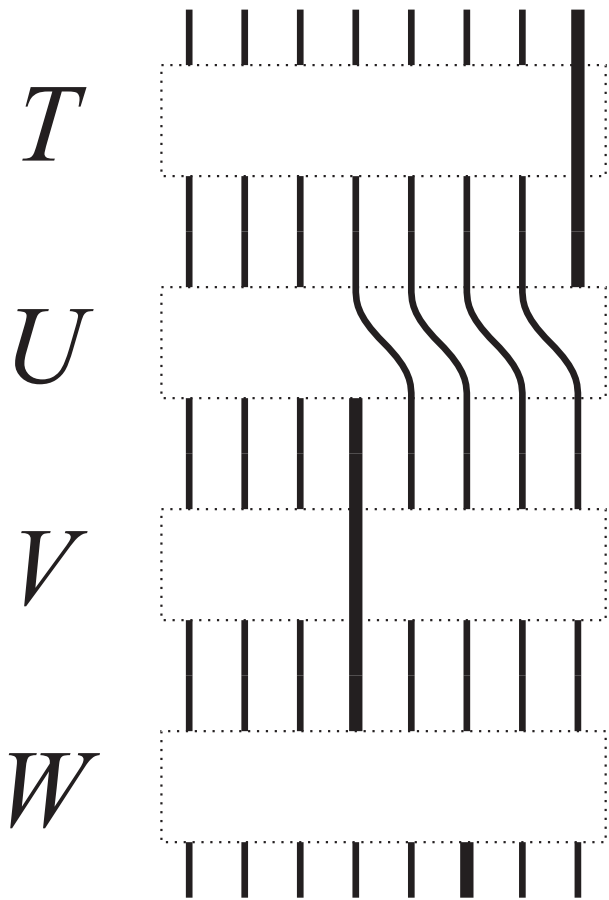} \caption{The $TUVW$-decomposition.} \label{pic:6}
\end{figure}

Now, an arbitrary braid word $w$ can be represented as follows:
\[
T \cdot U \cdot V \cdot W
\]
where the corresponding braid diagram takes the form shown on
Figure \ref{pic:6}. Note that between $T$ and $U$, $U$ and $V$,
$V$ and $W$ no free cancellation occurs. I.e. $T$  is a word in
which only small crossings occur. In $U$ every crossing involves
the $n+1$-th  strand and therefore is big. The word $V$ is a
letter, the corresponding crossing $\rho$ is small and $W$ is the
rest of the braid word $w$. We shall further move the only
crossing in $V$, that is $\rho$, upwards (towards the beginning of
the word). Thereby obtaining the following word $TV\cdot W'$,
where $TV$ is a word in which every crossing is small and $W'$
counts as many small crossings as $W$.

We prove the statement using double induction. First is on the
length $j$ of the word $U$ (the number of big crossings preceding
small crossing $\rho$) and the second one on the number $k$ of
letters in $W$ that correspond to small crossings (the number of
small crossings that stand further from the beginning of the word
than $\rho$).

Consider two cases.

1. The letter $V$ commutes with the terminal letter in $U$. In
which case we permute the letter $V$ and the terminal letter in
$U$ then freely reduce the word and obtain the new word $w=T\cdot
U'\cdot V\cdot U'' W$ (or $w=T'UW$, provided that $UV=VU$) and,
therefore, the number $j$ is decreased and the statement follows
by induction.

2. The terminal letter in $U$  does not commute with $V$. Since
every crossing in $T$ and $V$ is small and every crossing in $U$
is big the only possibilities for $V$ and the two crossings
preceding $V$ are shown on Figure \ref{pic:pv}. Replace the
configurations shown on Figure \ref{pic:pv} by the corresponding
braid on Figure \ref{pic:repl} and then freely reduce the obtained
word.

\newpage
Denote by $U'''$ the following word $U''' \cdot z_1 \cdot z_2=U$,
where $z_1, z_2 \in \left\{ x_1, \ldots, x_n, x_1^{-1}, \ldots,
x_n^{-1} \right\}$. In terms of words the transformations shown on
Figures \ref{pic:pv} and \ref{pic:repl} can be rewritten as
follows.
\begin{equation} \label{eq:I1} TU''' \cdot
x_{i+1}^{\delta}x_i^{\epsilon}x_{i+1}^{\epsilon} W \rightarrow
TU''' \cdot x_i^{\epsilon}x_{i+1}^{\epsilon}x_i^{\delta} W
\end{equation}
\begin{equation} \label{eq:I2}
TU''' \cdot x_{i+1}^{\epsilon}x_i^{\epsilon}x_{i+1}^{\delta} W
\rightarrow TU''' \cdot
x_i^{\delta}x_{i+1}^{\epsilon}x_i^{\epsilon} W
\end{equation}
\begin{equation} \label{eq:I3}
TU''' \cdot x_{i+1}^{\epsilon}x_{i+1}^{\epsilon}x_i^{\delta} W
\rightarrow TU''' \cdot
x_i^{\delta}x_{i+1}^{\delta}x_i^{\epsilon}x_i^{\epsilon}x_{i+1}^{-
\delta} W
\end{equation}
\begin{equation} \label{eq:I4}
TU''' \cdot x_{i+1}^{\epsilon}x_i^{-\epsilon}x_{i+1}^{\epsilon} W
\rightarrow TU''' \cdot
x_i^{\epsilon}x_{i+1}^{\epsilon}x_i^{\epsilon}x_i^{\epsilon}x_{i+1}^{-
\epsilon}x_{i+1}^{- \epsilon}x_i^{-\epsilon} W
\end{equation}
\begin{equation} \label{eq:II1}
TU''' \cdot x_i^{\delta}x_{i+1}^{\epsilon}x_i^{\epsilon} W
\rightarrow TU''' \cdot
x_{i+1}^{\epsilon}x_i^{\epsilon}x_{i+1}^{\delta} W
\end{equation}
\begin{equation} \label{eq:II2}
TU''' \cdot x_i^{\epsilon}x_{i+1}^{\epsilon}x_i^{\delta} W
\rightarrow TU''' \cdot
x_{i+1}^{\delta}x_i^{\epsilon}x_{i+1}^{\epsilon} W
\end{equation}
\begin{equation} \label{eq:II3}
TU''' \cdot x_i^{\epsilon}x_i^{\epsilon}x_{i+1}^{\delta} W
\rightarrow TU''' \cdot
x_{i+1}^{\delta}x_i^{\delta}x_{i+1}^{\epsilon}x_{i+1}^{\epsilon}x_i^{-
\delta} W
\end{equation}
\begin{equation} \label{eq:II4}
TU''' \cdot x_i^{\epsilon}x_{i+1}^{- \epsilon}x_i^{\epsilon} W
\rightarrow TU''' \cdot
x_{i+1}^{\epsilon}x_i^{\epsilon}x_{i+1}^{\epsilon}x_{i+1}^{\epsilon}x_i^{-
\epsilon}x_i^{- \epsilon}x_{i+1}^{-\epsilon} W.
\end{equation}

Furthermore, in the constructed procedure whenever the letter $V$
commutes with the terminal letter in $U$ we permute the letter $V$
and the terminal letter in $U$. In this case the word $U$ may
consist of a single letter. This transformation rewrites in terms
of words as follows:
\begin{equation} \label{eq:com1}
T\tilde{U}z_1V W \rightarrow T\tilde{U} V z_1   W,
\end{equation}
where $\tilde{U}z_1=U$ and $\tilde{U}$ is possibly empty.

Finally on the frontier of either $U'''$ (or $\tilde{U}$) or $W$
or $z_1$ and the replaced segment free cancellation might occur.
For the sake of completeness, we write the expression of free
reduction
\begin{equation} \label{eq:cancel}
x_i^{\epsilon} x_i^{-\epsilon} \rightarrow 1, \ \epsilon =\pm 1
\end{equation}
We leave the reader to check that Equations
(\ref{eq:I1})-(\ref{eq:II4}), (\ref{eq:com1}) and
(\ref{eq:cancel}) are in fact equalities in the braid group.

Suppose that after a step of the procedure from the initial word
$w$ we obtain a new word $w'$, and consider its $T\cdot U\cdot
V\cdot W$-decomposition. In the obtained $T\cdot U\cdot V\cdot
W$-decomposition the length $j$ of $U$ is lower, provided that the
crossing $\rho$ is not in $T$ (in which case the statement follows
by induction on $k$). The word $V$ is a letter and it corresponds
to the same crossing $\rho$, provided that the letter $V$ has not
cancelled with a letter from $T$. Finally, since all of the above
replacements do not create any new small crossings in $w$ (though
replacements  given by Expressions (\ref{eq:I3}), (\ref{eq:I4}),
(\ref{eq:II3}), (\ref{eq:II4}) create new big crossings), the
above process stops and the statement follows by induction on $k$.

As an output we obtain the following decomposition of an arbitrary
braid $w=\xi_n w_{n+1}$, where in the word $\xi_n$ every crossing
is small and in $w_{n+1}$ all big crossings of $w$ are collected.
We next apply a similar process to $\xi_n$ and take it in the form
in which all crossings involving the $n$-th strand are gathered in
its end, $w=\xi_n w_{n+1}=\xi_{n-1} w_n w_{n+1}$ and so on.
Finally, we obtain a word in the form (\ref{eq:NF}).

\subsection{Uniqueness}

Below we show that if $\xi_1$ and $\xi_2$ are two braid words
written in the form (\ref{eq:NF}) and $\xi_1=\xi_2$ in the braid
group $B_{n+1}$ then $\xi_1 \equiv \xi_2$. Firstly, note that, on
account of Theorem \ref{thm:art} (Theorem 17, \cite{Artin}) or
Theorem 6 in Section 11, \cite{Markov} if $\xi_1$ and $\xi_2$ are
two pure braid words written in the form (\ref{eq:NF}) and
$\xi_1=\xi_2$ then $\xi_1 \equiv \xi_2$.

Consider two braid words $\xi_1=p_2 \cdots p_{n}u$ and
$\xi_2=q_2q_3 \cdots q_{n}w$ written in the form (\ref{eq:NF}).
Denote $p_2 \cdots p_{n}$ and $q_2q_3 \cdots q_{n}$ by,
correspondingly, $t$ and $v$. Note that $t$ and $v$ consist of
small crossings only, while $u$ and $w$ of big crossings only.

In the above notation we have
\begin{equation} \label{eq:uniq}
t\circ u \circ w^{-1} \circ v^{-1}= 1 \hbox{ in what follows that
} v^{-1}\circ t\circ u \circ w^{-1}=1,
\end{equation}
where by `$\circ$' we denote the concatenation (without free
reduction) of words.

We next want to show that the subword $v^{-1}\circ t$ of the word
$v^{-1}\circ t\circ u\circ w^{-1}$ is a pure braid every crossing
in which is small and that the subword $u\circ w^{-1}$ of the word
$v^{-1}\circ t\circ u\circ w^{-1}$ is a pure braid every crossing
in which is big. Obviously $v^{-1}t$ is a braid in which every
crossing is small (by the definition both $v$ and $t$ involve
letters $x_1, \ldots, x_{n-1}$ and their inverse only).

\begin{remark} \label{rem:sb}
Let $\xi$ be a braid every crossing in which is big {\rm(}in which
case the initial letter of $\xi$ is $x_n${\rm)} and let $\zeta$ be
an arbitrary braid every crossing in which is small {\rm(}in which
case $\zeta$ does not contain the letter $x_n^{\pm 1}${\rm)}. Then
in the word $\zeta \cdot \xi$ every crossing in the subword
$\zeta$ is small, while every crossing in the subword $\xi$ is
big.
\begin{figure}[h]
\centering
  \includegraphics[height=2in, width=12.6cm, keepaspectratio]{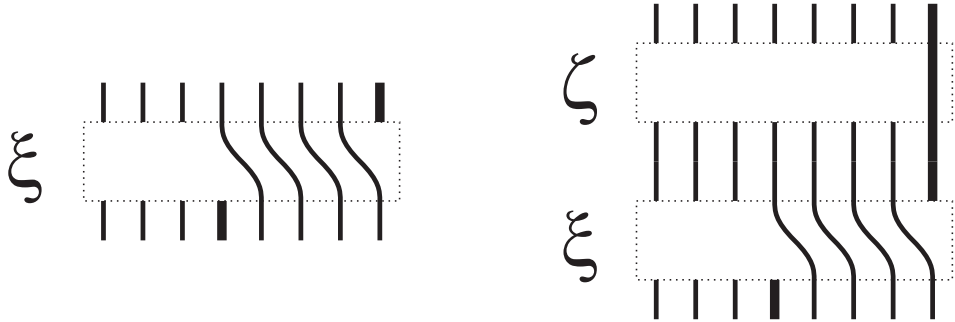}
  \caption{The braid $\xi$ is big regardless of the initial segment of the word
  $\zeta \xi$, provided that every crossing in $\zeta$ is small}
    \label{pic:rem}
\end{figure}
By the definition the braid $\xi$ entangles the $n+1$-th strand,
leaving the first $n$ strands parallel {\rm(}regardless of the
permutation defined by the subword $\zeta$ of the word $\zeta
\xi${\rm)}, see Figure {\rm \ref{pic:rem}}.
\end{remark}

We next show that every crossing in the subword $u\circ w^{-1}$ of
the word $v^{-1}\circ t\circ u\circ w^{-1}$ is big. Since every
crossing in $u$ is big, by Remark \ref{rem:sb}, every crossing in
the subword $u$ of the word $v^{-1}\circ t \circ u$ is big. We are
left to show that every crossing in the $w^{-1}$ passage of the
word $v^{-1}\circ t\circ u\circ w^{-1}$ is big. To show the latter
we need to study the structure of crossings of the inverse of a
braid that entangles the $n+1$-th strand only. Consider a braid
$\xi$ so that every its crossing is big (it therefore entangles
the $n+1$-th strand, leaving the first $n$ strands unentangled)
and so that the permutation defined by $\xi$ takes the $n+1$-th
strand into the $k$-th position then the inverse $\xi^{-1}$
entangles the $k$-th strand only, leaving the other $n$ strands
unentangled. Furthermore, by the definition, the permutation
defined by $\xi^{-1}$ is the inverse of the permutation defined by
$\xi$ and consequently the permutation defined by $\xi^{-1}$ takes
the $k$-th strand of $\xi^{-1}$ into the $n$-th position.

\begin{figure}[h]
\centering
  \includegraphics[height=3cm, width=6cm, keepaspectratio]{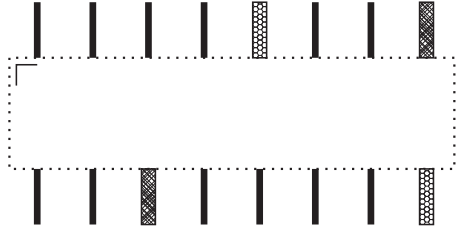} \caption{An arbitrary braid}
    \label{pic:w}
\end{figure}

\begin{figure}[h]
\centering
  \includegraphics[height=3cm, width=6cm, keepaspectratio]{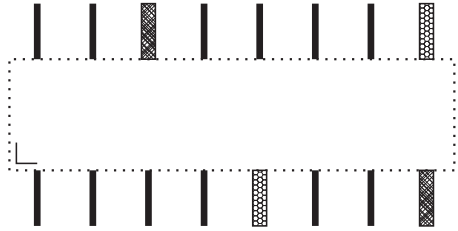} \caption{The inverse of an arbitrary braid}
    \label{pic:winv}
\end{figure}

\begin{figure}[h]
\centering
  \includegraphics[height=1.23in, width=0.75in, keepaspectratio]{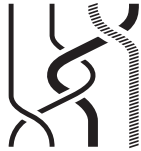} \caption{The diagram of $\xi = x_3x_2^{-2}x_1$}
    \label{pic:expl2}
\end{figure}

\begin{figure}[h]
\centering
  \includegraphics[height=1.19in, width=0.81in, keepaspectratio]{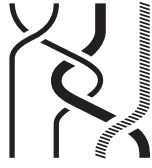}
  \caption{The diagram of $\xi^{-1} =
x_1^{-1}x_2^{2}x_3^{-1}$}
    \label{pic:expl3a}
\end{figure}

The above argument can also be illustrated by the following
diagrams. The diagram of an arbitrary braid word $\xi$ has the
form given on Figure \ref{pic:w} (two strands are designated and a
corner is marked to demonstrate the connections between a braid
and its inverse)
and, consequently, the diagram of the word $w^{-1}$ has the form
given on Figure \ref{pic:winv},
i.e. the diagram of the word $\xi^{-1}$ is a two times reflected
diagram of $\xi$. Let us consider an example.
\begin{example}
Let $\xi = x_3x_2^{-2}x_1$, the respective diagram is given by
Figure {\rm\ref{pic:expl2}}. The inverse of $\xi$ is $\xi^{-1} =
x_1^{-1}x_2^{2}x_3^{-1}$ and the corresponding diagram is given by
Figure {\rm\ref{pic:expl3a}}.
\end{example}

We now return to the consideration of the word $v^{-1}\circ t\circ
u\circ w^{-1}$. Recall that, by our assumption, the word $u$ takes
the $n+1$-th strand into the $k$-th position. Since $tu=vw$ thus
the permutations defined by $tu$ and $vw$ coincide, and both $u$
and $w$ permute the $n+1$-th strand into the $k$-th position. Now,
from the foregoing discussion we know that:
\begin{itemize}
    \item Every crossing in $v^{-1} \circ t$ is small,
    \item Every crossing in $u$ is big,
    \item The word $w^{-1}$ entangles the $k$-th strand only.
\end{itemize}
Consequently, in the word $v^{-1}\circ t\circ u\circ w^{-1}$ the
subword $w^{-1}$ entangles the strand, which is on the $k$-th
position. By Remark \ref{rem:sb}, since every crossing in the
braid $v^{-1}t$ is small and since the permutation defined by the
word $u$ takes the $n+1$-th strand into the $k$-th position, so
does the permutation defined by the word $v^{-1}\circ t\circ u$.
We thereby obtain that all crossings in the $v^{-1}\circ t$
subword of the word $v^{-1}\circ t\circ u\circ w^{-1}$ are small
and all the crossings in the subword $u\circ w^{-1}$ are big.

\begin{figure}[h]
\centering
  \includegraphics[height=4in, width=12.6cm, keepaspectratio]{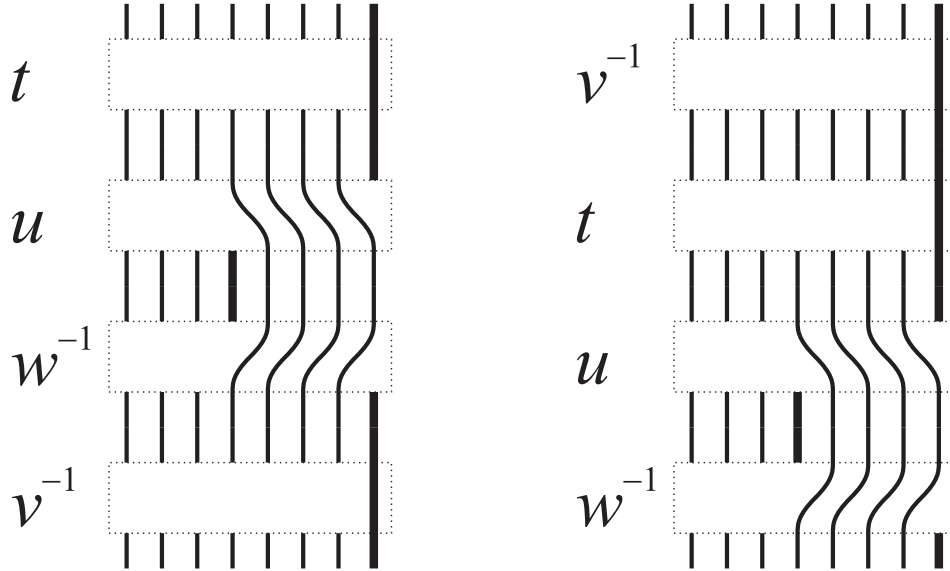}
  \caption{Diagrams of $t\circ u\circ w^{-1}\circ v^{-1}$ and $t\circ u \circ w^{-1}\circ v^{-1}$}
    \label{pic:diag}
\end{figure}

The latter argument can be illustrated by diagrams. The diagrams
corresponding to Equation (\ref{eq:uniq}) take the form shown on
Figure \ref{pic:diag}.
From Figure \ref{pic:diag} we see that, since the permutation
defined by the word $v^{-1}\circ t\circ u\circ w^{-1}$ is the
identity and since $v^{-1}\circ t$ defines a permutation of the
first $n$ strands, while $u\circ w^{-1}$ entangles the $n+1$-th
strand only both $v^{-1}\circ t$ and $u\circ w^{-1}$ are pure
braids.

From the above discussion, since $v^{-1}\circ t\circ u\circ
w^{-1}$ is a pure braid, by Theorem \ref{thm:art} (Theorem 17,
\cite{Artin}) or Theorem 6 in Section 11, \cite{Markov}, the
equality $v^{-1}\circ t\circ u\circ w^{-1}=1$ therefore implies
that
\[
v^{-1}\circ t=1 \hbox{ and } u\circ w^{-1}=1.
\]
Using induction on the number of strands $n+1$ we may assume that
the equality $v=t$ in $B_n$ implies that $v \equiv t$.

Finally, we are left to show that if the braid $uw^{-1}$ is
trivial, every its crossing is big and the word $uw^{-1}$ is
freely reduced then $uw^{-1}$ is the empty word. Since $uw^{-1}$
is a pure braid every crossing of which is big, it can be written
as a product of the following words
\[
A_{i,n+1}=x_n^{-1} \cdots x_{i+1}^{-1} \cdot x_{i}^2 \cdot
x_{i+1}\cdots x_n;\; i=1,\dots ,n.
\]
By the definition $A_{n,n+1}=x_n^2$. Consequently normal form
(\ref{eq:art}) of a pure braid $uw^{-1}$ is $U_1$ (in the notation
of Theorem \ref{thm:art}). Therefore, since $uw^{-1}=1$, on behalf
of Theorem \ref{thm:art} (Theorem 17, \cite{Artin}) or Theorem 6
in Section 11, \cite{Markov}, $uw^{-1}$ is the empty word and $u
\equiv w$. Consequently, normal form (\ref{eq:NF}) is unique.

\section{Actual Implementation and Random Braid} \label{sec:rb}

In this section we explain another approach to the computation of
normal form (\ref{eq:NF}). This approach turns out to be efficient
in the actual implementation of the algorithm on a computer. The
current and the following sections are less formal and hold some
of our ideas and comments on the issue.

In the actual implementation of the algorithm for computing the
normal form (\ref{eq:NF}) of a braid word $\xi$ it is more
convenient to treat elements of $B_{n+1}$ as sequences of
crossings. Let us consider an example.

\begin{example}
Let $\xi = x_3x_2^{-2}x_1$. This element rewrites in terms of
crossings as follows:
\[
\xi_1 = \mid 3,4 \mid \mid 2,4 \mid^{-1} \mid 2,4 \mid^{-1} \mid
1,2 \mid.
\]
On the other hand, given a proper sequence of crossings, we can
form the corresponding braid word. Let us consider the following
sequence of crossings:
\[
\xi_1^*=\mid 1, 2 \mid \mid 3,4\mid \mid 1,4 \mid \mid 2,4 \mid
^{-1} \mid 2,4 \mid ^{-1} \mid 1,4 \mid ^{-1}.
\]
It is fairly obvious that the corresponding braid word is $\xi ^*=
x_1x_3x_2x_1^{-2}x_2^{-1}$.
\end{example}

\begin{remark}
The correspondence between sequences of crossings and braid words
is not one-to-one. Let $\zeta_1= \mid 1,3 \mid ^{-1}$ be a
sequence of crossings. There is no braid word $\zeta$
corresponding to the sequence $\zeta_1$.
\end{remark}

\begin{remark} \label{rem:seqcrreg}
The set of all sequences of crossings that represent elements of a
fixed braid group $B_n$ is regular, i. e. recoginsed by a finite
automaton.
\end{remark}
\begin{proof}
The automaton has $n!$ states which are numbered by the elements
of the group of permutations $S_n$ on $n$ symbols. The initial
state is the one that corresponds to the trivial permutation, and
all the states are fail states. There are two edges, which are
labelled by crossings $\mid i, j \mid ^{\pm 1}$  from a state
labelled $\theta\in S_n$ to a state labelled $\vartheta\in S_n$
whenever $\vartheta = \theta \cdot (i,j)=\theta \cdot
{(l,l+1)}^\theta=(l,l+1)\theta$, for some $l=1,\ldots, n$. We
leave the reader to check that this automaton recognises the set
of all sequences of crossings that correspond to braid words (not
necessarily freely reduced).
\end{proof}

We next rewrite the replacements given by Figures \ref{pic:pv} and
\ref{pic:repl} (Equations (\ref{eq:I1}) - (\ref{eq:II4})) in terms
of sequences of crossings.

Let $k>j,l$, $1\le k,j,l\le n+1$ and let $\epsilon, \delta \in
\left\{ -1, 1 \right\}$. Replacements given by Figures
\ref{pic:pv} and \ref{pic:repl} (Equations (\ref{eq:I1}) -
(\ref{eq:II4})) take the form (we assume that $k>j,l$):
\begin{equation}\label{eq:I1'}
\mid j,k\mid^{\delta}\mid l,k\mid^{\epsilon}\mid
l,j\mid^{\epsilon}
 \rightarrow
\mid l,j\mid^{\epsilon}\mid l,k\mid^{\epsilon}\mid
j,k\mid^{\delta}
\end{equation}
\begin{equation}\label{eq:I2'}
\mid j,k\mid^{\epsilon}\mid l,k\mid^{\epsilon}\mid
l,j\mid^{\delta}
 \rightarrow
\mid l,j\mid^{\delta}\mid l,k\mid^{\epsilon}\mid
j,k\mid^{\epsilon}
\end{equation}
\begin{equation}\label{eq:I3'}
\mid j,k\mid^{\epsilon}\mid j,k\mid^{\epsilon}\mid
l,j\mid^{\delta} \rightarrow \mid l,j\mid^{\delta}\mid
l,k\mid^{\delta}\mid j,k\mid^{\epsilon} \mid
j,k\mid^{\epsilon}\mid l,k\mid^{-\delta}
\end{equation}
\begin{equation}\label{eq:I4'}
\mid j,k\mid^{\epsilon}\mid l,k\mid^{-\epsilon}\mid
l,j\mid^{\epsilon} \rightarrow \mid l,j\mid^{\epsilon}\mid
l,k\mid^{\epsilon}\mid j,k\mid^{\epsilon}\mid
j,k\mid^{\epsilon}\mid l,k\mid^{-\epsilon} \mid
l,k\mid^{-\epsilon}\mid j,k\mid^{-\epsilon}
\end{equation}

Recall that in Subsection \ref{ss:exist} we used yet another
transformation of a braid. In the $T\cdot U \cdot V \cdot W$
decomposition of a braid if the letter $V$ commutes with the
terminal letter in $U$ we permute the letter $V$ and the terminal
letter in $U$ and obtain the new word $w=T U'VU'' W$, see Equation
(\ref{eq:com1}). To perform such a transformation we introduce the
following transformation of sequences of crossings.

Let $i<\max \left\{l,k\right\}$ and $j<\max\left\{l,k\right\}$;
$1\le i,j,k,l \le n+1$ and let $\epsilon, \delta \in \left\{ -1, 1
\right\}$. The transformation mentioned in Subsection
\ref{ss:exist} takes the form:
\begin{equation} \label{eq:com2}
\mid l,k \mid^{\epsilon}\mid i,j \mid^{\delta}\rightarrow \mid l,k
\mid^{\epsilon}\mid i,j \mid^{\delta}
\end{equation}
We, therefore replace the sequence of crossings $\rho\varrho$ by
$\varrho\rho$, whenever $\rho$ and $\varrho$ commute and $\varrho$
is `smaller' than $\rho$.

Finally, since in Subsection \ref{ss:exist} we used reduction in
the free group (see Equation (\ref{eq:cancel})) we need to
introduce the following transformation of sequences of crossings:
\begin{equation} \label{eq:D}
\mid i,j \mid^{\epsilon} \mid i,j \mid^{- \epsilon} \rightarrow 1.
\end{equation}

We next reformulate Theorem \ref{thm:NF} for elements of braid
groups viewed as sequences of crossings.
\begin{theorem} \label{thm:NFcr}
Every sequence of crossings which corresponds to an element $w \in
B_{n+1}$ can be taken {\rm(}by the means of  transformations
{\rm(\ref{eq:I1'})--(\ref{eq:D}))} to the form
\begin{equation} \label{eq:NFcr}
w=\mid 1,2\mid^m \cdot  w_3(\mid 1,3\mid,\mid 2, 3\mid) \cdots
w_{n+1}(\mid 1,n\mid, \dots ,\mid n,n+1\mid);\ m \in \mathbb{Z}
\end{equation}
Where for every $3 \le k \le n+1$ the words $w_k$ do not contain
subwords of the form $\mid i,j\mid^{\epsilon} \mid
i,j\mid^{-\epsilon}$, $\epsilon=\pm 1$. Under the above
assumptions presentation in the form {\rm(\ref{eq:NFcr})} is
unique.
\end{theorem}

\begin{remark} \label{rem:preserve}
Let $w$ be an arbitrary sequence of crossings, which represents a
braid word $\xi$. Suppose next that in $\xi$ there is a subword,
which coincides with one of the left-hand sides of rules
{\rm(\ref{eq:I1'})--(\ref{eq:D})}. Then the word $\xi'$, which is
obtained from $\xi$ by replacing the left-hand side of a rule
{\rm(\ref{eq:I1'})--(\ref{eq:D})} by the respective right-hand
side, is a sequence of crossings which rewrites into a braid word.
I. e. rules {\rm(\ref{eq:I1'})--(\ref{eq:D})} preserve the
property of being a representative of a braid word.
\end{remark}

\begin{remark}
Not only did the above approach turned out to be very useful in
actual implementation of the algorithm of computation of normal
form {\rm(\ref{eq:NF})}, but also hints at us the idea to
construct a term rewriting system {\rm(}a Knuth-Bendix like
algorithm{\rm)} for elements of braid groups viewed as sequences
of crossings. We suppose that  normal form {\rm(\ref{eq:NF})} can
be generalised and a similar normal form can be constructed for an
arbitrary Artin group $A$ of finite type. First one needs to fix
an ordering of the natural generators of the corresponding
{\rm(}to $A${\rm)} Coxeter group $\mathcal{A}$. Next one needs to
rewrite an arbitrary word $\zeta \in A$ as a sequence of elements
of $\mathcal{A}$. Then transform the element $\zeta$
{\rm(}agreeing with the introduced order{\rm)} to a form whose
presentation as a sequence of generators of $\mathcal{A}$ is in
some sense small {\rm(}an analogue of the form
{\rm(\ref{eq:NF}))}. And consequently there exist a nice geometric
normal form for an arbitrary Artin group of finite type. Below, in
Section {\rm\ref{sec:ex}}, we show how one can construct such a
form for the group $A=\left<a,b \mid abab=baba\right>$.
\end{remark}

One of the most important problems in a struggle to construct
cryptography on braid groups is to give a method for generating a
random braid (see \cite{DehSurv}). We, therefore, can not but
notice that Theorem \ref{thm:NF} gives a method of generating a
random braid.

A naive approach to generating a random braid is to generate a
random freely reduced word and claim this a random braid. However
there is a dramatic difference between a random braid and a random
freely reduced word. For example, the results of \cite{DeN} and
\cite{Vershik} show that a random braid has got a nontrivial
centraliser (its centraliser differs from the center of $B_{n+1}$)
with a non-zero probability, while computational results show,
that the centraliser of a random reduced word in the free group,
treated as an element of $B_{n+1}$ is trivial.

To generate a random braid on $n+1$ strands we can use Theorem
\ref{thm:NF}, for instance, as follows. Here we use the ideas of
\cite{mult}, whereto we refer the reader for details and
justification of the method. We do not in any way insist that the
suggested method is better than any other method known, though we
believe that it may turn out to be useful for computer scientists.

We suggest the following random process with weak interferences.
One has:
\[
w=x_1^m \cdot w_3(x_1,x_2) \cdot w_4(x_1,x_2,x_3) \cdots
w_{n+1}(x_1, \dots ,x_n).
\]

\textbf{1}. We generate a random power of $x_1$ as follows, thus
entangling the first two strands. We start at the identity element
of $B_{n+1}$ and either do nothing with probability $s_2\in (0,1]$
and then go to step \textbf{2} or move to one of the two elements
$x_1$ or $x_1^{-1}$ with equal probabilities $\frac{1-s_2}{2}$. If
we are at an element $x_1^m=v\ne 1$ we either stop at $v$ with
probability $s_2$ (and proceed to step \textbf{2}), or move, with
probability $1-s_2$ to the vertice $x_1^{m+1}$, if $m>0$ and to
the vertice $x_1^{m-1}$, if $m<0$.

\textbf{2}. On this step we generate $w_3(x_1,x_2)$ (in the
notation of Equation (\ref{eq:NF})), thus entangling the third
strand into the first two. This process can be treated as a
discrete random walk on two points, which are linked to the first
and the second strands. We start at the identity element of
$B_{n+1}$ and either do nothing with probability $s_3\in (0,1]$
and then go to step \textbf{3} or move to one of the two elements
$x_2$ or $x_2^{-1}$ with equal probabilities $\frac{1-s_3}{2}$. If
we are at an element $x_2^{m_1}\cdot x_1^{m_2} \cdots x_l^{m_t}=
v\ne 1$, $l=1,2$ we either stop at $v$ with probability $s_3$ (and
proceed to step \textbf{3}), or:
\newline
\textbf{A}. If $l=2$ and $m_t>0$ is odd, move, with probability
$\frac{1-s_3}{3}$, to one of the vertices $x_2^{m_1}\cdot
x_1^{m_2} \cdots x_2^{m_t+1}$ or $x_2^{m_1}\cdot x_1^{m_2} \cdots
x_l^{m_t}x_1$ or $x_2^{m_1}\cdot x_1^{m_2} \cdots
x_l^{m_t}x_1^{-1}$.
\newline
\textbf{B}. If $l=2$ and $m_t<0$ is odd, move, with probability
$\frac{1-s_3}{3}$, to one of the vertices $x_2^{m_1}\cdot
x_1^{m_2} \cdots x_2^{m_t-1}$ or $x_2^{m_1}\cdot x_1^{m_2} \cdots
x_l^{m_t}x_1$ or $x_2^{m_1}\cdot x_1^{m_2} \cdots
x_l^{m_t}x_1^{-1}$.
\newline
\textbf{C}. If $l=2$ and $m_t>0$ is even, move, with probability
$1-s_3$, to the vertice $x_2^{m_1}\cdot x_1^{m_2} \cdots
x_2^{m_t+1}$.
\newline
\textbf{D}. If $l=2$ and $m_t<0$ is even,  move, with probability
$1-s_3$, to the vertice $x_2^{m_1}\cdot x_1^{m_2} \cdots
x_1^{m_t-1}$.
\newline
\textbf{E}. If $l=1$ and $m_t>0$ is odd,  move, with probability
$1-s_3$, to the vertice $x_2^{m_1}\cdot x_1^{m_2} \cdots
x_1^{m_t+1}$.
\newline
\textbf{F}. If $l=1$ and $m_t<0$ is odd,  move, with probability
$1-s_3$, to the vertice $x_2^{m_1}\cdot x_1^{m_2} \cdots
x_1^{m_t-1}$.
\newline
\textbf{G}. If $l=1$ and $m_t<0$ is even,  move, with probability
$\frac{1-s_3}{3}$, to one of the vertices $x_2^{m_1}\cdot
x_1^{m_2} \cdots x_1^{m_t-1}$ or $x_2^{m_1}\cdot x_1^{m_2} \cdots
x_1^{m_t}x_2$ or $x_2^{m_1}\cdot x_1^{m_2} \cdots x_1^{m_t}
x_2^{-1}$.
\newline
\textbf{H}. If $l=1$ and $m_t>0$ is even,  move, with probability
$\frac{1-s_3}{3}$, to one of the vertices $x_2^{m_1}\cdot
x_1^{m_2} \cdots x_1^{m_t+1}$  or $x_2^{m_1}\cdot x_1^{m_2} \cdots
x_1^{m_t}x_2$ or $x_2^{m_1}\cdot x_1^{m_2} \cdots x_1^{m_t}
x_2^{-1}$.

\textbf{3}. On this step we generate $w_4(x_1,x_2,x_3)$ (in the
notation of Equation (\ref{eq:NF})), thus entangling the fourth
strand into the first three. This process can be treated as a
discrete random walk on three points, which are linked to the
first three strands. We start at the identity element of $B_{n+1}$
and either do nothing with probability $s_4\in (0,1]$ and then go
to step \textbf{4} or move to one of the two elements $x_3$ or
$x_3^{-1}$ with equal probabilities $\frac{1-s_4}{2}$. If we are
at an element $v_1 \cdot x_l^{m}= v\ne 1$ we either stop at $v$
with probability $s_4$ (and proceed to step \textbf{4}), or:
\newline
\textbf{A}. If $l=3$ and $m>0$ is odd, move, with probability
$\frac{1-s_3}{3}$, to one of the vertices $v_1 \cdot x_3^{m+1}$ or
$v_1 \cdot x_3^{m}x_2$ or $v_1 \cdot x_3^{m}x_2^{-1}$.
\newline
\textbf{B}. If $l=3$ and $m<0$ is odd, move, with probability
$\frac{1-s_3}{3}$, to one of the vertices $v_1 \cdot x_3^{m-1}$ or
$v_1 \cdot x_3^{m}x_2$ or $v_1 \cdot x_3^{m}x_2^{-1}$.
\newline
\textbf{C}. If $l=3$ and $m>0$ is even, move, with probability
$1-s_3$, to the vertices $v_1 \cdot x_3^{m+1}$.
\newline
\textbf{D}. If $l=3$ and $m<0$ is even, move, with probability
$1-s_3$, to the vertices $v_1 \cdot x_3^{m-1}$.
\newline
\textbf{E}. If $l=2$ and $m>0$ is odd, move, with probability
$\frac{1-s_3}{3}$, to one of the vertices $v_1 \cdot x_2^{m+1}$ or
$v_1 \cdot x_2^{m}x_1$ or $v_1 \cdot x_2^{m}x_1^{-1}$.
\newline
\textbf{F}. If $l=2$ and $m<0$ is odd, move, with probability
$\frac{1-s_3}{3}$, to one of the vertices $v_1 \cdot x_2^{m-1}$ or
$v_1 \cdot x_2^{m}x_1$ or $v_1 \cdot x_2^{m}x_1^{-1}$.
\newline
\textbf{G}. If $l=2$ and $m>0$ is even, move, with probability
$\frac{1-s_3}{3}$, to one of the vertices $v_1 \cdot x_2^{m+1}$ or
$v_1 \cdot x_2^{m}x_3$ or $v_1 \cdot x_2^{m}x_3^{-1}$.
\newline
\textbf{H}. If $l=2$ and $m<0$ is even, move, with probability
$\frac{1-s_3}{3}$, to one of the vertices $v_1 \cdot x_2^{m-1}$ or
$v_1 \cdot x_2^{m}x_3$ or $v_1 \cdot x_2^{m}x_3^{-1}$.
\newline
\textbf{I}. If $l=1$ and $m>0$ is odd, move, with probability
$1-s_3$, to the vertice $v_1 \cdot x_1^{m+1}$.
\newline
\textbf{J}. If $l=1$ and $m<0$ is odd, move, with probability
$1-s_3$, to the vertice $v_1 \cdot x_1^{m-1}$.
\newline
\textbf{K}. If $l=1$ and $m>0$ is even, move, with probability
$\frac{1-s_3}{3}$, to one of the vertices $v_1 \cdot x_1^{m+1}$ or
$v_1 \cdot x_1^{m}x_2$ or $v_1 \cdot x_1^{m}x_2^{-1}$.
\newline
\textbf{L}. If $l=1$ and $m<0$ is even, move, with probability
$\frac{1-s_3}{3}$, to one of the vertices $v_1 \cdot x_1^{m-1}$ or
$v_1 \cdot x_1^{m}x_2$ or $v_1 \cdot x_1^{m}x_2^{-1}$.
\newline

\textbf{4}. And so on.

In other words on the \textbf{$k$}-th step we generated the word
$w_{k+1}(x_1,\ldots,x_k)$, thus entangling the $k+1$-th strand
into the first $k$. Since the words $w_j$'s from Equation
(\ref{eq:NF}) are determined uniquely and since, by the
construction, the obtained words $w_j$'s are freely reduced, we
obtain a `random' braid.

\section{A Term Rewriting System for Elements of Braid Groups}

In this section we view elements of $B_{n+1}$ as sequences of
crossings, which we for the most part will refer as strings. We
construct a confluent string rewriting system, whose output will
be a string in the form (\ref{eq:NFcr}).

To construct a rewriting system one is to define the strings and
the rules. In our case, strings are the sequences of crossings
that correspond to braid words and the rules of the rewriting
system are the rules given by (\ref{eq:I1'})--(\ref{eq:D}). We
also term these rules by transformations or rewrites.

Let us number all the crossings for the group $B_{n+1}$ as
follows: $\rho_1=\mid 1,2 \mid, \rho_2=\mid 1,3 \mid, \rho_3=\mid
2,3 \mid, \rho_4=\mid 1,4 \mid, \ldots,
\rho_{\frac{n(n+1)}{2}}=\mid n, n+1 \mid$; naturally
$\rho_1^{-1}=\mid 1,2 \mid^{-1},\ldots,
\rho_{\frac{n(n+1)}{2}}^{-1}=\mid n, n+1 \mid^{-1}$.

We now prove that the rewriting system defined on the set of all
sequences of crossings that correspond to braids by
transformations (\ref{eq:I1'})--(\ref{eq:D}) is confluent, i. e.
that any chain of rewrites terminates and that there is the only
residue and that it is the one corresponds to normal form
(\ref{eq:NFcr}).

 To prove the confluence of our rewriting system we first prove
that for any string one can sequentially apply only a finite
number of transformations. First of all note that all rewrites
(\ref{eq:I1'})--(\ref{eq:D}) either have the form
\begin{equation} \label{eq:rules}
g_1  (\rho_{k_j}, \dots , \rho_{k_m})\cdot
\rho_{k_i}^{\epsilon}\rightarrow \rho_{k_i}^{\epsilon} \cdot g_2
(\rho_{k_j}, \dots , \rho_{k_m}),
\end{equation}
where $g_1 (\rho_{k_j}, \dots , \rho_{k_m}),g_2 (\rho_{k_j}, \dots
, \rho_{k_m})$ are strings of crossings that involve crossings
$\rho_{k_j},\ldots, \rho_{k_m}$ and their inverses only, $k_i<k_j,
\ldots, k_j$, or the form
\begin{equation} \label{eq:rules2}
\rho_i^\epsilon \rho_i^{-\epsilon} \rightarrow 1,
\end{equation}
where $i=1,\ldots, \frac{n(n+1)}{2}$ and $\epsilon =\pm 1$.

Denote by $n(w)$ the maximal length of a chain of transformations
applied to a string $w$, so that $n(w)$ is either a positive
integer or the symbol $\infty$, in the case that the length of
chains of transformations is not bounded above.

\begin{remark}
Let $w=v\cdot u$ be a presentation of a sequence of crossings $w$
so that no free cancellation between $v$ and $u$ occurs and
$n(w)\ne \infty$. Then $n(w)\ge n(v),n(u)$.
\end{remark}

We now intend to prove the following

\begin{lemma} \label{lem:bound}
Let $R$ be the set of rules of the form {\rm(\ref{eq:rules})} and
{\rm(\ref{eq:rules2})}. Let $w$ be an arbitrary string then any
sequence of transformations of the string $w$ terminates.
\end{lemma}
\begin{proof}
Suppose first that a string involves $\rho_{\frac{n(n+1)}{2}}$ and
its inverse only. Then any chain of transformations has
transformations given by Equation {\rm(\ref{eq:rules2})} only, and
is clearly finite. Thus, without loss of generality, we may
assume that any chain of transformations of any sequence of
crossings that involves $\rho_2,\ldots,\rho_{\frac{n(n+1)}{2}}$
(and their inverses) only counts a finite number of steps.

We now use induction on the number of occurrences of $\rho_1^{\pm
1}$ in a string. If $\rho_1^{\pm 1}$ does not occur in a string,
the statement is straightforward. Consider next a sequence of
crossings
$t=w(\rho_2,\ldots,\rho_m)\rho_1^{\epsilon}v(\rho_2,\ldots,\rho_m)$,
$m=\frac{n(n+1)}{2}$. By the induction assumption, $n(w)$ and
$n(v)$ are finite. Consider an arbitrary chain of  transformations
of $t$. Suppose that the $t$-th rewrite involves the
$\rho_1^{\epsilon_1}$ and the string rewrites as follows
$w'(\rho_2,\ldots,\rho_m)\rho_1^{\epsilon}v'(\rho_2,\ldots,\rho_m)\rightarrow
w''(\rho_2,\ldots,\rho_m)\rho_1^{\epsilon}v''(\rho_2,\ldots,\rho_m)$.
Where $w''(\rho_2,\ldots,\rho_m)$ is a substring of
$w'(\rho_2,\ldots,\rho_m)$ and $n(w'') \le n(w')<\infty$ and
$n(v'')<\infty$. Applying a similar argument to the word
$w''(\rho_2,\ldots,\rho_m)\rho_1^{\epsilon}v''(\rho_2,\ldots,\rho_m)$
and noticing that $n(w'')<n(w)$ or the length of $w''$ is strictly
lower than the one of $w$ (for example if $w=w'$) or both, we
conclude that any chain of transformations of this sequence of
crossings terminates.

Consider next an arbitrary string with $k$ occurrences of
$\rho_1$ and write it in the form:
\[
w=w_1(\rho_2,\ldots,\rho_m)\rho_1^{\epsilon_1}w_2(\rho_2,\ldots,\rho_m)\rho_1^{\epsilon_2}\cdots
w_k(\rho_2,\ldots,\rho_m)\rho_1^{\epsilon_k}w_{k+1}(\rho_2,\ldots,\rho_m).
\]
By induction we may assume that for the sequence of crossings
$n(\bar
w=w_1(\rho_2,\ldots,\rho_m)\rho_1^{\epsilon_1}w_2(\rho_2,\ldots,\rho_m)\rho_1^{\epsilon_2}\cdots
w_k(\rho_2,\ldots,\rho_m))$ is finite. Consider an arbitrary chain
of transformations of $w$. Suppose that its $t$-th step involves
the $\rho_1^{\epsilon_k}$ and the string rewrites as follows
$w'(\rho_1,\rho_2,\ldots,\rho_m)\rho_1^{\epsilon_k}w'_{k+1}(\rho_2,\ldots,\rho_m)\rightarrow
w''(\rho_1,\rho_2,\ldots,\rho_m)
\rho_1^{\epsilon_k}w''_{k+1}(\rho_2,\ldots,\rho_m)$, where $w''$
is a subword of $w'$ and thus $n(w'') \le n(w')<\infty$, and
$n(w''_{k+1}(\rho_2,\ldots,\rho_m))<\infty$. Applying a similar
argument to  $w''\rho_1^{\epsilon_k}w_{k+1}(\rho_2,\ldots,\rho_m)$
the statement follows. Indeed, either $n(w'')$ is strictly lower
than $n(w')$ or the length of $w''$ is strictly lower than the
length of $\bar w$ or both and $n(w''_{k+1})$ is always finite,
the procedure, therefore, terminates.
\end{proof}

\begin{theorem}
The set of rules {\rm(\ref{eq:I1'})--(\ref{eq:D})} gives rise to a
confluent term rewriting system on the set of sequences of
crossings, which represent braid words.
\end{theorem}
\proof By the definition, a rewriting system is confluent whenever
the order of application of rules does not matter, i. e. we have
the following diagrams:
\[
\begin{array}{ccccccc}
                                         &            &  x\cdot r_1\cdot y \cdot l_2 \cdot z &          & \\
                                         & \nearrow  &                                       & \searrow &  \\
w=x\cdot \underbrace{l_1}_{\hbox{a rule}} \cdot y \cdot \overbrace{l_2}^{\hbox{a rule}} \cdot  z  &  & &  & \texttt{res}(w);\\
                                         & \searrow   &                                      & \nearrow & \\
                                         &            &  x\cdot l_1\cdot y \cdot r_2 \cdot z &          & \\
\end{array}
\]
\[
\begin{array}{ccccccc}
                                         &            &  x\cdot r_1\cdot l'_2 \cdot z &          & \\
                                         & \nearrow  &                                       & \searrow &  \\
w=x \cdot \lefteqn{\underbrace{  \phantom{l'_1 \cdot y}}_{\hbox{a rule}}} l'_1\cdot \overbrace{ y\cdot l'_2}^{\hbox{a rule}} \cdot z & & & & \texttt{res}(w);\\
                                         & \searrow   &                                      & \nearrow & \\
                                         &            &  x\cdot l'_1\cdot r_2 \cdot z &          & \\
\end{array}
\]
\[
\begin{array}{ccccccc}
                                         &            &  x\cdot r_1 \cdot z & & \\
                                         & \nearrow  &                                       & \searrow &  \\
w=x\cdot \underbrace{l'_1 \cdot \overbrace{l_2}^{\hbox{a rule}} \cdot l''_1}_{\hbox{a rule}} \cdot  z  & & & & \texttt{res}(w).\\
                                         & \searrow   &                                      & \nearrow & \\
                                         &            &  x\cdot l'_1\cdot r_2 \cdot l''_1 \cdot z & & \\
\end{array}
\]
Here $l_1$ and $l_2$, $l'_1 \cdot y$ and $y\cdot l'_2$, $l'_1
\cdot l_2 \cdot l''_1$ and $l_2$ are left parts of the rules,
$r_1$ and $r_2$ are right-hand sides of the rules and
$\verb"res"(w)$ is the residue of $w$, i.e. such a word to which
no transformation can be applied. This is roughly saying that if
there is a choice of which rule to apply then there exists a
common residue of the resulting words.

To show that our term rewriting system is confluent suppose that
there exist two distinct residues $\texttt{res}_1(w)$ and
$\texttt{res}_2(w)$ of a string $w$. By Theorem \ref{thm:NFcr}
there exist a residue $w^*$ of $w$ in the form (\ref{eq:NFcr})
(since none  of the rules (\ref{eq:I1'})--(\ref{eq:D}) can be
applied to $w^*$). Since $\texttt{res}_1(w)\ne \texttt{res}_2(w)$
at least one of them does not coincide with $w^*$. Suppose that
$\texttt{res}_1(w)\ne w^*$. Then, by Remark \ref{rem:preserve},
$\texttt{res}_1(w)$ corresponds to a braid, and therefore, on
account of Theorem \ref{thm:NF} (see also Section \ref{ss:exist}),
we can take $\texttt{res}_1(w)$ to $w^*$ by the means of the rules
(\ref{eq:I1'})--(\ref{eq:D})
--- a contradiction. \hfill $\blacksquare$

\begin{corollary}
The set of all sequences of crossings that correspond to normal
form (\ref{eq:NFcr}) of elements of the braid group $B_{n+1}$ is
regular.
\end{corollary}
\begin{proof}
By Remark \ref{rem:seqcrreg} the set of all sequences of crossings
that correspond to elements of $B_{n+1}$ is regular. Since the
collection of rules {\rm(\ref{eq:I1'})--(\ref{eq:D})} is finite,
the set of all sequences of crossings that do not contain the
left-hand sides of rules is regular. The statement now follows,
for the intersection of two regular sets is regular.
\end{proof}

\section{An Analogue of the Normal Form for another Artin Group}
\label{sec:ex}

In this Section we show how one can obtain an analogue of the form
(\ref{eq:NF}) for an Artin group of finite type $A=\left<a,b \mid
abab=baba\right>$.

Any Artin group of type $B$ embeds into a braid group on
sufficiently many strands and the embedding is fairly natural, see
\cite{Vershinin} for details. For example, the group $A$ embeds
into $B_3$ as follows: $\phi(a)=x_1, \phi(b)=x_2^2$. Using this
embedding and the normal form constructed in this paper one can
obtain similar geometric normal forms for $A$ (and any other Artin
group of type $B$). However, we want to explicitly demonstrate the
generalisation of the idea used to construct normal forms in braid
groups and show how one can construct normal forms in $A$
`barehanded'.

Consider the corresponding Coxeter group $\mathcal{A}$. As it is
well-known, see for instance \cite{Bourb}, $\mathcal{A}$ can be
treated as a group, generated by the following reflections of a
square: $\bar a=\mid 1,3 \mid$ and $\bar b=\mid 1,2\mid \mid
3,4\mid$, where the corresponding square has the form:
\[
\begin{CD}
1   @=  2   \\
@|      @|  \\
4   @=  3
\end{CD}
\]
We use this notation to express, that $\bar a$ is the reflection
of the square that permutes vertices $1$ and $3$ and $\bar b$ is
the reflection of the square that permutes $1$ with $2$, and $3$
with $4$:
\[
a:
\begin{CD}
1   @=  2   \\
@|      @|  \\
4   @=  3
\end{CD}
\longrightarrow
\begin{CD}
3   @=  2   \\
@|      @|  \\
4   @=  1
\end{CD}
\hbox{ and } b:
\begin{CD}
1   @=  2   \\
@|      @|  \\
4   @=  3
\end{CD}
\longrightarrow
\begin{CD}
2   @=  1   \\
@|      @|  \\
3   @=  4
\end{CD}
\]

We next consider every element of $A$ as a sequence of elements of
$\mathcal{A}$. The group $\mathcal{A}$ consists of $8$ elements:
$1,\mid 1,3 \mid$, $\mid 2,4 \mid$, $\mid 1,2 \mid \mid 3,4 \mid$,
$\mid 1,3 \mid \mid 1,2 \mid \mid 3,4 \mid$, $\mid 1,2 \mid \mid
3,4 \mid  \mid 1,3 \mid$, $\mid 1,4 \mid \mid 2,3 \mid$, $\mid 1,3
\mid \mid 2,4 \mid$. All the above elements, but $\mid 1,3\mid$
involve the $4$-th vertice of the square.

\begin{remark}
In the case of braid groups, the crossing $\mid 1,3 \mid$ of two
strands can correspond to any generator {\rm(}since all the
permutations $(i,j)$ are conjugate in the Coxeter group, which
corresponds to the braid group, i. e. the group of permutations on
$n+1$ symbols{\rm)}. However, only the elements $\mid 1,3 \mid$
and $\mid 2,4 \mid$ of the group $\mathcal{A}$ can correspond to
the generator $a$ {\rm(}treated as a subword of a word of
$A${\rm)}, since $\mid 1,3 \mid$ is the image of $a$ in
$\mathcal{A}$ and $\mid 2,4 \mid$ is the only element of
$\mathcal{A}$ conjugate to $\mid 1,3 \mid$.
\end{remark}

We next prove the following theorem
\begin{theorem}
In the above notation, every element $w \in A$ can be taken to the
form
\begin{equation} \label{eq:nf1}
w=a^m \cdot w_1(a,b); m \in \mathbb{Z}
\end{equation}
Where $w_1$ is freely reduced and the reflections, which
correspond to the letters in $w_1$ are the following $\mid 2,4
\mid$, $\mid 1,2 \mid \mid 3,4 \mid$, $\mid 1,3 \mid \mid 1,2 \mid
\mid 3,4 \mid$, $\mid 1,2 \mid \mid 3,4 \mid 1,3 \mid$, $\mid 1,4
\mid \mid 2,3 \mid$, $\mid 1,3 \mid \mid 2,4 \mid$ only.
\end{theorem}
\proof The proof of the theorem repeats the argument given in
Subsection \ref{ss:exist}. We only need to introduce analogues of
transformations (\ref{eq:I1})-(\ref{eq:II4}) and Figures
\ref{pic:pv} and \ref{pic:repl}. Algebraic analogues of Figures
\ref{pic:pv} and \ref{pic:repl} take the form:
\begin{gather} \notag
\begin{split}
& (\mid 1,2 \mid \mid 3,4 \mid)^{\epsilon}\cdot(\mid 1,2 \mid \mid
3,4 \mid)^{\epsilon} \cdot \mid 1,3 \mid ^{\delta} \rightarrow
\mid 1,3 \mid ^{\delta} \cdot (\mid 1,4  \mid \mid 2,3
\mid)^{\delta} \cdot
\\
\cdot \mid 2,4 \mid^{\delta} \cdot & (\mid 1,2 \mid \mid 3,4
\mid)^{\epsilon} \cdot (\mid 1,2 \mid \mid 3,4 \mid)^{\epsilon}
\cdot \mid 2,4 \mid^{-\delta} \cdot (\mid 1,2\mid \mid 3,4
\mid)^{-\delta};
\end{split}
\end{gather}
\begin{gather} \notag
\begin{split}
& (\mid 1,4 \mid \mid 2,3 \mid)^{\epsilon} \cdot(\mid 1,4 \mid
\mid 2,3 \mid)^{\epsilon} \cdot \mid 1,3 \mid ^{\delta}
\rightarrow
 \mid 1,3 \mid^{\delta} \cdot (\mid 1,2  \mid \mid 3,4
\mid)^{\delta} \cdot \\ \cdot \mid 2,4 \mid^{\delta} \cdot & (\mid
1,4 \mid \mid 2,3 \mid)^{\epsilon} \cdot (\mid 1,4 \mid \mid 2,3
\mid)^{\epsilon} \cdot\mid 2,4 \mid^{-\delta} \cdot (\mid 1,4 \mid
\mid 2,3 \mid)^{-\delta};
\end{split}
\end{gather}
\begin{gather} \notag
\begin{split}
\mid 2,4 \mid^{\epsilon}&\cdot (\mid 1,4 \mid \mid 2,3
\mid)^{\epsilon} \cdot \mid 1,3 \mid^{\epsilon} \rightarrow\\
\rightarrow &(\mid 1,2 \mid \mid 3,4 \mid)^{\epsilon} \cdot \mid
1,3 \mid^{\epsilon} \cdot (\mid 1,4 \mid \mid 2,3
\mid)^{\epsilon}\cdot \mid 2,4 \mid^{\epsilon} \cdot (\mid 1,2
\mid \mid 3,4 \mid)^{-\epsilon};
\end{split}
\end{gather}
\begin{gather} \notag
\begin{split}
\mid 2,4 \mid^{\epsilon}&\cdot (\mid 1,2 \mid \mid 3,4
\mid)^{\epsilon} \cdot \mid 1,3 \mid^{\epsilon} \rightarrow\\
\rightarrow & (\mid 1,4 \mid \mid 2,3 \mid)^{\epsilon} \cdot \mid
1,3 \mid^{\epsilon} \cdot (\mid 1,2 \mid \mid 3,4
\mid)^{\epsilon}\cdot \mid 2,4 \mid^{\epsilon} \cdot (\mid 1,4
\mid \mid 2,3 \mid)^{-\epsilon};
\end{split}
\end{gather}
\begin{gather} \notag
\begin{split}
\mid 2,4 \mid^{\epsilon}&\cdot (\mid 1,4 \mid \mid 2,3
\mid)^{-\epsilon} \cdot \mid 1,3 \mid^{-\epsilon} \rightarrow\\
\rightarrow (\mid 1,2 \mid \mid 3,4 \mid)^{-\epsilon}&  \cdot \mid
1,3 \mid^{-\epsilon} \cdot (\mid 1,4 \mid \mid 2,3
\mid)^{-\epsilon}\cdot \mid 2,4 \mid^{\epsilon} \cdot (\mid 1,2
\mid \mid 3,4 \mid)^{\epsilon};
\end{split}
\end{gather}
\begin{gather} \notag
\begin{split}
\mid 2,4 \mid^{\epsilon}&\cdot (\mid 1,2 \mid \mid 3,4
\mid)^{-\epsilon} \cdot \mid 1,3 \mid^{-\epsilon} \rightarrow\\
\rightarrow (\mid 1,4 \mid \mid 2,3 \mid)^{-\epsilon}& \cdot \mid
1,3 \mid^{-\epsilon} \cdot (\mid 1,2 \mid \mid 3,4
\mid)^{-\epsilon}\cdot \mid 2,4 \mid^{\epsilon} \cdot (\mid 1,4
\mid \mid 2,3 \mid)^{\epsilon};
\end{split}
\end{gather}
\begin{gather} \notag
\begin{split}
\mid 2,4 \mid^{-\epsilon}&\cdot (\mid 1,4 \mid \mid 2,3
\mid)^{-\epsilon} \cdot \mid 1,3 \mid^{\epsilon} \rightarrow\\
\rightarrow(\mid 1,2 \mid \mid 3,4 \mid)^{\epsilon}& \cdot \mid
1,3 \mid^{\epsilon} \cdot (\mid 1,4 \mid \mid 2,3
\mid)^{-\epsilon}\cdot \mid 2,4 \mid^{-\epsilon} \cdot (\mid 1,2
\mid \mid 3,4 \mid)^{-\epsilon};
\end{split}
\end{gather}
\begin{gather} \notag
\begin{split}
\mid 2,4 \mid^{-\epsilon}&\cdot (\mid 1,2 \mid \mid 3,4
\mid)^{-\epsilon} \cdot \mid 1,3 \mid^{\epsilon} \rightarrow\\
\rightarrow(\mid 1,4 \mid \mid 2,3 \mid)^{\epsilon}& \cdot \mid
1,3 \mid^{\epsilon} \cdot (\mid 1,2 \mid \mid 3,4
\mid)^{-\epsilon}\cdot \mid 2,4 \mid^{-\epsilon} \cdot (\mid 1,4
\mid \mid 2,3 \mid)^{-\epsilon};
\end{split}
\end{gather}
\begin{gather} \notag
\begin{split}
\mid 2,4 \mid^{\epsilon}&\cdot (\mid 1,4 \mid \mid 2,3
\mid)^{-\epsilon} \cdot \mid 1,3 \mid^{\epsilon} \rightarrow (\mid
1,2 \mid \mid 3,4 \mid)^{\epsilon} \cdot \mid 1,3 \mid^{\epsilon}
\cdot (\mid 1,4 \mid \mid 2,3 \mid)^{\epsilon}\cdot\\
\cdot\mid 2,4 \mid^{\epsilon} &\cdot \mid 2,4 \mid^{\epsilon}
\cdot (\mid 1,4 \mid \mid 2,3 \mid)^{-\epsilon}\cdot  (\mid 1,4
\mid \mid 2,3 \mid)^{-\epsilon}\cdot\mid 2,4 \mid^{-\epsilon}
\cdot (\mid 1,2 \mid \mid 3,4 \mid)^{-\epsilon};
\end{split}
\end{gather}
\begin{gather} \notag
\begin{split}
\mid 2,4 \mid^{\epsilon}&\cdot (\mid 1,2 \mid \mid 3,4
\mid)^{-\epsilon} \cdot \mid 1,3 \mid^{\epsilon} \rightarrow
 (\mid 1,4 \mid \mid 2,3 \mid)^{\epsilon} \cdot \mid 1,3 \mid^{\epsilon} \cdot (\mid 1,2 \mid \mid 3,4
\mid)^{\epsilon}\cdot
\\
\cdot\mid 2,4 \mid^{\epsilon}&\cdot \mid 2,4 \mid^{\epsilon}
\cdot(\mid 1,2 \mid \mid 3,4 \mid)^{-\epsilon} \cdot (\mid 1,2
\mid \mid 3,4 \mid)^{-\epsilon}\cdot\mid 2,4 \mid^{-\epsilon}\cdot
(\mid 1,4 \mid \mid 2,3 \mid)^{-\epsilon}.
\end{split}
\end{gather}
In terms of words analogues of transformations
(\ref{eq:I1})-(\ref{eq:II4}) take the form:
\begin{equation}  TU''' \cdot
b^{\epsilon}b^{\epsilon}a^{\delta} W \rightarrow TU''' \cdot
a^{\delta}b^{\delta}a^{\delta}  b^{\epsilon}b^{\epsilon}
a^{-\delta}a^{-\delta} W
\end{equation}
\begin{equation}
TU''' \cdot a^{\epsilon}b^{-\epsilon}a^{\epsilon} W \rightarrow
TU''' \cdot
b^{\epsilon}a^{\epsilon}b^{\epsilon}a^{\epsilon}a^{\epsilon}
b^{-\epsilon}b^{-\epsilon}a^{-\epsilon}b^{-\epsilon} W
\end{equation}
\begin{equation}
TU''' \cdot a^{\epsilon}b^{\epsilon}a^{\epsilon} W \rightarrow
TU''' \cdot
b^{\epsilon}a^{\epsilon}b^{\epsilon}a^{\epsilon}b^{-\epsilon} W
\end{equation}
\begin{equation}
TU''' \cdot a^{\epsilon}b^{-\epsilon}a^{-\epsilon} W \rightarrow
TU''' \cdot
b^{-\epsilon}a^{-\epsilon}b^{-\epsilon}a^{\epsilon}b^{\epsilon} W
\end{equation}
\begin{equation}
TU''' \cdot a^{-\epsilon}b^{-\epsilon}a^{\epsilon} W \rightarrow
TU''' \cdot
b^{\epsilon}a^{\epsilon}b^{-\epsilon}a^{-\epsilon}b^{-\epsilon} W
\end{equation}
One can see that each rule in the list above moves the reflection
$\mid 1, 3 \mid^{\pm 1}$ to the beginning of the word and is so
that in the right-hand side of the rules there is only one
occurrence of $a$, which corresponds to the reflection $\mid 1, 3
\mid^{\pm 1}$. Furthermore, the rules above enumerate all the
possibilities, when $\mid 1, 3 \mid^{\pm 1}$  corresponds to the
terminal letter of a subword of a word from $A$ and so that only
the following elements of $\mathcal{A}$ correspond to the $2$
letters preceding the letter that corresponds to $\mid 1,3
\mid^{\pm 1}$:
\[\mid 2,4 \mid, \mid 1,2 \mid \mid 3,4 \mid, \mid
1,3 \mid \mid 1,2 \mid \mid 3,4 \mid, \mid 1,2 \mid \mid 3,4 \mid
1,3 \mid, \mid 1,4 \mid \mid 2,3 \mid, \mid 1,3 \mid \mid 2,4
\mid. \] \hfill $\blacksquare$

\begin{remark}
We did not show that the form {\rm(\ref{eq:nf1})} is unique. To
prove the uniqueness of normal form {\rm(\ref{eq:NF})} for braids
we used results of E. Artin. We do not know similar results for
arbitrary Artin groups of finite type. However, it may be derived
from the fact that $A$ embeds into some braid group.
\end{remark}

\section{Questions and Final Remarks}

In the current section we list some of the questions that arose in
our work. We also briefly discuss some of advantages and
disadvantages of normal form (\ref{eq:NF}).

\begin{question}
In Section {\rm\ref{sec:rb}} we have given a description of a
procedure of generating a random braid. That procedure involves
$n$ random walks, each of which is defined by a parameter $s_j$.
Choose the $s_j$'s so that the obtained word generator would be
useful for practical needs of computer scientists. We suppose that
a solution of this problem can be piloted by papers
{\rm\cite{DeN}} and {\rm\cite{Vershik}}.
\end{question}

Consider a list of transformations $( \blacklozenge )$ of the
words of a free monoid with the alphabet $X \cup X^{-1}$.
\begin{itemize}
    \item $x_i^{\epsilon}x_i^{-\epsilon} \rightarrow
    1, \; \epsilon=\pm 1;$
    \item $x_i^{\epsilon}x_j^{\eta} \rightarrow
    x_j^{\eta}x_i^{\epsilon}, \; |i-j| \ge 2; \epsilon, \eta =
    \pm 1;$
    \item $x_i^{\epsilon}x_{i+1}^{\epsilon}x_i^{\epsilon}
    \leftrightarrow
    x_{i+1}^{\epsilon}x_i^{\epsilon}x_{i+1}^{\epsilon}, \; \epsilon=
    \pm 1;$
    \item $x_i^{\epsilon}x_{i+1}^{k \epsilon}x_i^{\eta}
    \leftrightarrow
    x_{i+1}^{\eta}x_i^{k \epsilon}x_{i+1}^{\epsilon}, \; \epsilon, \eta =
    \pm 1, \epsilon \cdot \eta = -1, k \in \mathbb{N};$
    \item $x_{i+1}^{\epsilon}x_i^{k \epsilon}x_{i+1}^{\eta}
    \leftrightarrow
    x_i^{\eta}x_{i+1}^{k \epsilon}x_i^{\epsilon}, \; \epsilon, \eta =
    \pm 1, \epsilon \cdot \eta = -1, k \in \mathbb{N};$
\end{itemize}

\begin{question} \label{q:2}
I. Consider two geodesic words $\xi_1,\xi_2$ from a braid group
$B_n$ such that $\xi_1\stackrel{B_n}{=}\xi_2$. Is the list of
 transformations $( \blacklozenge )$ sufficient to  obtain the word
$\xi_1$ from the word $\xi_2$?

II. Let $\xi_1,\xi_2$ be arbitrary words from a braid group $B_n$.
Assume that $\xi_1$ is written in geodesic form and
$\xi_1\stackrel{B_n}{=}\xi_2$. Is the list of
 transformations $( \blacklozenge )$ sufficient to  obtain the word
$\xi_1$ from the word $\xi_2$?

III. Analogues of questions I and II for an arbitrary Artin group.
\end{question}

A positive answer to parts I and II of Question \ref{q:2} was
announced in \cite{Makanin}. However, as far as the authors are
concerned, proof never appeared.

\begin{question} \label{q:3}
We presume that the introduced normal forms do not form an
automatic structure {\rm(}this was conjectured in our conversation
with G. A. Noskov{\rm)}. Let $\xi$ be a braid word. Our
experiments show that normal form {\rm(\ref{eq:NF})} $\xi^*$ of
$\xi$ may have an exponential length on the length of the input
word {\rm(}consider the word $\xi=x_3^2 x_2^2 x_1^2 x_2^2$ and its
powers, for instance{\rm)}. In what follows that this normal form
can not be automatic {\rm(}see {\rm\cite{Epstein})}. However,
resulting from our experiments, we conjecture that in the class of
words conjugate to $\xi$ {\rm(}and even among all cyclic
permutation of $\xi${\rm)} there exists an element $\zeta$ so that
its normal form $\zeta^*$ is polynomial {\rm(}linear{\rm)} on the
length of its geodesic form $\zeta_1$.
\end{question}

\begin{question}
As mentioned in Question {\rm \ref{q:3}} normal form {\rm
(\ref{eq:NF})} $\xi^*$ of $\xi$ may have an exponential length on
the length of the input word. Is it true that normal form
{\rm(\ref{eq:NF})} of $\xi$ is generically polynomial
{\rm(}linear{\rm)} on the length of a geodesic word $\xi_1=\xi$?
For a detailed explanation of the term `generically polynomial'
the reader may consult {\rm\cite{kapovich}}.
\end{question}

In \cite{Artin} E. Artin wrote the following regarding his normal
forms. We can not but agree.
\begin{quote}
Though it has been proved that every braid can be transformed
 to a similar normal form the writer is convinced that any attempt to
carry this out on a living person would only lead to violent
protests and discrimination against mathematicians. He would
therefore discourage such an experiment. \flushright \emph{E.
Artin}
\end{quote}

The authors would like to thank Vladimir Vershinin and Referee of
the paper for their useful comments that helped to substantially
improve the text.

\end{document}